\documentclass[reqno]{jgokova}
\usepackage{graphicx}
\usepackage{graphics}
\usepackage{amsmath}
\usepackage{amscd}
\usepackage{latexsym}

\newcommand{\ben}{\begin{enumerate}}
\newcommand{\een}{\end{enumerate}}
\newcommand{\be}{\begin{equation}}
\newcommand{\ee}{\end{equation}}
\newcommand{\bea}{\begin{eqnarray}}
\newcommand{\eea}{\end{eqnarray}}
\newcommand{\bc}{\begin{center}}
\newcommand{\ec}{\end{center}}

\theoremstyle{definition}

\theoremstyle{remark}

 \volume{3}

\begin{document}

\newtheorem{Thm}{Theorem}
\newtheorem{Lem}[Thm]{Lemma}
\newtheorem{Cor}[Thm]{Corollary}
\newtheorem{Prop}[Thm]{Proposition}
\newtheorem{Rm}{Remark}

\def\a{{\mathbb a}}
\def\C{{\mathbb C}}
\def\A{{\mathbb A}}
\def\B{{\mathbb B}}
\def\D{{\mathbb D}}
\def\E{{\mathbb E}}
\def\R{{\mathbb R}}
\def\P{{\mathbb P}}
\def\S{{\mathbb S}}
\def\Z{{\mathbb Z}}
\def\O{{\mathbb O}}
\def\H{{\mathbb H}}
\def\V{{\mathbb V}}
\def\Q{{\mathbb Q}}
\def\Cn{${\mathcal C}_n$}
\def\CM{\mathcal M}
\def\CG{\mathcal G}
\def\CH{\mathcal H}
\def\CT{\mathcal T}
\def\CF{\mathcal F}
\def\CA{\mathcal A}
\def\CB{\mathcal B}
\def\CD{\mathcal D}
\def\CP{\mathcal P}
\def\CS{\mathcal S}
\def\CZ{\mathcal Z}
\def\CE{\mathcal E}
\def\CL{\mathcal L}
\def\CV{\mathcal V}
\def\CW{\mathcal W}
\def\IC{\mathbb C}
\def\IF{\mathbb F}
\def\IK{\mathcal K}
\def\IL{\mathcal L}
\def\IP{\bf P}
\def\IR{\mathbb R}
\def\IZ{\mathbb Z}

\title{An infinite family of exotic Dolgachev surfaces without $1$- and $3$- handles }
\author[AKBULUT]{Selman Akbulut}
\thanks{The author is partially supported by NSF grant DMS 0905917}
\address{Department  of Mathematics, Michigan State University,  MI, 48824}
\email{akbulut@math.msu.edu }
\date{\today}
\begin{abstract} 
Starting with the Dolgachev surface  $E(1)_{2,3}$ we construct an infinite \linebreak
family of 
exotic  copies of the rational surface $E(1)$, each of which admits a handlebody decomposition without $1$- and $3$- handles, and  we draw these handlebodies. 
\end{abstract}

\date{}
\maketitle

\setcounter{section}{-1}


\section{Introduction}

It is an old problem to determine when an exotic copy of a smooth simply connected \linebreak
$4$-manifold admits a handle decomposition without $1$- and $3$- handles. Note that if an exotic $S^4$ or $\C\P^2$ exists then its handle decomposition must contain either $1$- or \linebreak 
$3$-handles.  Finding such exotic manifolds realizing the smallest Betti number is a particularly 
 interesting problem. For example in  \cite{y} and also \cite{ay} exotic manifolds without \linebreak
  $1$- and \mbox{$3$-handles} were demonstrated. An interesting difficult case has been the Dolgachev \linebreak 
 surface $E(1)_{2,3}$, which is an exotic copy of $E(1)={\bf \C\P}^2 \#\; 9 \overline{{\bf \C\P}}^2$. In  \cite{a1} a handlebody of this manifold without $1$- and $3$-handles were drawn (previously this manifold was conjectured not to admit such a handlebody \cite{hkk}). Here by extending the technique of \cite{a1} we construct infinitely many exotic copies of the rational surface $E(1)$ without $1$- and $3$-handles and draw their handlebodies. 

{\Thm There is an infinite family  $X_n$, $n=1,2, \dots$ of mutually non-diffeomorphic exotic copies  of $E(1)$ without $1$- and $3$-handles, such that  $X_{1}=E(1)_{2,3}$ (their handlebody pictures are given in Figure 18).}

\vspace{.05in}

Here is a brief outline of this paper:  A common technique to generate infinitely many exotic copies of a smooth $4$-manifold  $X$ is the ``Knot surgery" operation \mbox{$X\leadsto X_{K}$} of \cite{fs}, where $K\subset S^3$ is a knot.  Here the main goal is to make this operation compatible with the handle cancelling techniques of \cite{a1}.  At first try this seems to be an impossible task. The main difficulty with removing $1$- and $3$-handles from $E(1)_{K}$ is due to a certain twisting on its handles, which can only be seen after turning its handlebody upside down. To understand this better, we first start  with the test case of the handlebody of $E(1)_{K_{1}}=E(1)_{2,3}$, where $K_{1}$ is the trefoil knot, given in Figure 41 of \cite{a1}. By turning  this handlebody upside down we get a nice symmetric picture of it in Figure 15 (recall, the Figure 41 of \cite{a1} was already obtained from an initial handlebody picture of  $E(1)_{K_{1}}$ by turning it upside down twice!).  By inspecting the boundary of this handlebody  in Figure 16 we locate the twisting among the attaching circles of its handles, and in return this helps us to understand how this twist was compensated in the  handles of $E(1)_{K_{1}}$ so that no $1$- and $3$- handles are needed.  Then by imitating this  handle configuration in the more general case  in Figure 18, we can construct handlebodies for $E(K)_{K_{n}}$ so that no  $1$- and $3$-handles are needed, where $K_n$,  $n=1,2,..$ are knots with distinct Alexander polynomials (this implies these manifolds are mutually distinct).


\section{Construction}

We start with the first picture of Figure 1 which is the  handlebody of $E(1)_{K_{1}}$ given in  Figure 41 of \cite{a1}. By a handle slide we get the second picture of the Figure 1. Next we turn this handlebody upside down. This is done by finding a diffeomorphism from the boundary of this handlebody to $\partial B^4$ and by attaching the dual $2$-handles (the small \mbox{$0$-framed} blue linking circles in the first picture of Figure 2) to $B^4$ via this diffeomorphism. During this diffemorphism (e.g. sliding red handles over each other and blowing up and down operations) no red handles can slide over the dual blue handles, but the blue handles can slide over the red handles and they can slide over each other. 

\vspace{.05in}

Now we first replace a copy of $S^1\times B^3$, in the interior of the first picture of Figure 2, with $B^2\times S^2$ (i.e. we replace the ``dotted'' circle  with $0$-framed circle), then blow down the six $-1$ framed red circles to obtain the second picture of Figure 2 (this process increases the framing of the newly introduced $0$-framed circle to 6). Next by the indicated handle slides, isotopies and blowing downs we obtain Figures 3,4,5, $\dots$ and finally the first picture of Figure 13. Then by canceling a pair of $1$- and $2$-handles from this picture we obtain the second picture of Figure 13, and then by isotopies and the indicated handle slides we get Figure 14, and then end up with the first picture of Figure 15, which is another handlebody representation of $E(1)_{K_{1}}$. A careful inspection shows that this  handlebody has no $1$- and $3$-handles. It has no $3$-handles because its boundary is $S^3$ (this will be checked in the next paragraph), and  it has no $1$-handles because this figure has eleven $2$-handles and one ``slice'' $1$-handle (= two standard $1$-handles plus one $2$-handle), and two  $1$-handles of the slice $1$-handle are cancelled by the two $2$-handles corresponding to two of the small $1$-framed linking circles (this is exactly the same  canceling process performed in the first paragraph of the Section 3 of \cite{a1}). Hence this handlebody can be considered as a handlebody consisting of just ten $2$-handles.
 
 \vspace{.05in}

We can now check the boundary of the the first picture of Figure 15 is $S^3$ as follows:  First  we replace the slice $1$-handle with a $2$-handle (i.e. replace dot with zero) and slide it over the two $1$-framed two handles (corresponding to two trefoil knots) as indicated by the arrows. This gives the first picture of Figure 16, then by further indicated handle slides in Figure 16 and cancellations we get $S^3$. Note that this process creates a left twist to the handlebody which undoes the right twist on the upper most $2$- handle of the Figure~16. 

\vspace{.05in}

Having seen how the boundary of the first picture of Figure 15 is identified by $S^3$, we can  check that the $1$-framed dotted green circle (denoted by $\lambda$ in the picture) on the boundary is just a $0$-framed unknot. Hence we can add a $2$-handle to $\lambda$ with $1$-framing, and immediately after that cancel it with a $3$-handle. Now by  sliding over the handle $\lambda$ we get the second picture of the Figure 15, which is another more symmetric handlebody representation of $E(1)_{K_{1}}$ with eleven $2$-handles and one $3$-handle (after canceling the slice $1$-handle with $2$-handles). The advantage of this picture is that, it can easily be identified as the manifold obtained by doing the ``knot surgery operation"  to the first picture in Figure 17  by using the algorithm of \cite{a3} and \cite{a4}. Also the reader can easily check that the  first picture of Figure 17 is just E(1).

\vspace{.05in}

Now this knot surgery operation process $E(1)\leadsto E(1)_{K_{1}}$ can be imitated to construct $E(1)\leadsto E(1)_{K_{n}}$ with the sequence of knots $K_n$ (drawn in the second picture of Figure~17) with the same conclusions; i.e. the  first picture of Figure 18 is the knot surgered \linebreak
manifold $E(1)_{K_{n}}$, and the second picture of Figure 18 is a handlebody representation of this manifold  $E(1)_{K_{n}}$ without $1$- and $3$- handles. Let $X_{n}=E(1)_{K_{n}}$. \qed


{\Rm The Reader might wonder, why  once we found the  handlebody configuration of E(1) in Figure 17 we didn't shorten this paper by making it as the starting point of the proof (hence making the paper independent from \cite{a1}). There are two reasons: First we don't know if the cusp in Figure 17 is the cusp of the elliptic fibration of $E(1)$ (even though it can be checked that as homology classes they coincide). Hence  starting the proof with Figure 17  would give the conclusion of the Theorem 1  but not the identification of $X_{1}$ with $E_{1,2}$. So the proof of the first part of the Theorem 1 is independent of \cite{a1}. The second reason relating Theorem 1  to \cite{a1} is that we wanted to set a dictionary between the handles of $E(1)_{K_{n}}$ and $E(1)_{2,3}$, so in future we may be able to relate the corks in $E(1)_{K_{n}}$ with the cork  of $E(1)_{2,3}$ (\cite{ay}).  Recall that in \cite{a1} it was shown that  $E(1)_{2,3}$ is obtained from  $E(1)$ by twisting along the certain cork $\overline{W}_1\subset E(1)$. It is an interesting problem to determine the cork structures of the rest of the manifolds $X_n$, $n\in {\Z}.$}

\vspace{.05in}

{\Rm Notice  that $K_{-3}$ is a slice knot, i.e. $K_{-3} =\partial D_{0}^2$ for some disk $D_{0}^2\subset B^4$. Hence 
$K_{-3 } \;\# -K_{-3 }$ bounds two slice disks in $B^4$, one is $D^2_{1}= D_{0}^2\;\natural -D_{0}^2 \subset B^4$, and  the other one is the usual slice disk $D_{2}^2\subset B^4$  induced from the imbedding $ K_{-3 } \times I \subset B^2\times B^2 = B^4$. \linebreak
 Recall ``the slice $1$-handle'' is $B^4-D_{2}^2$. It can be checked that we replace the slice $1$-handle in Figure 18  by the other slice complement $B^4-D_{1}^2$ amazingly the resulting 
  manifold splits ${\C\P}^2 $'s, i.e. $X_n$ changes its smooth structure. This is related to the
  phenomena discussed in \cite{a5} (and also in \cite{a4}): A knot in $S^3$ can bound two different ribbon disks in $B^4$  such that the two ribbon complements are exotic copies of each other relative to their boundaries (and they can be taken to be homotopy equivalent to $S^1$).  The example given here may be viewed as the closed manifold version of this phenomenon.}

\begin{figure}
\includegraphics[width=.70\textwidth]{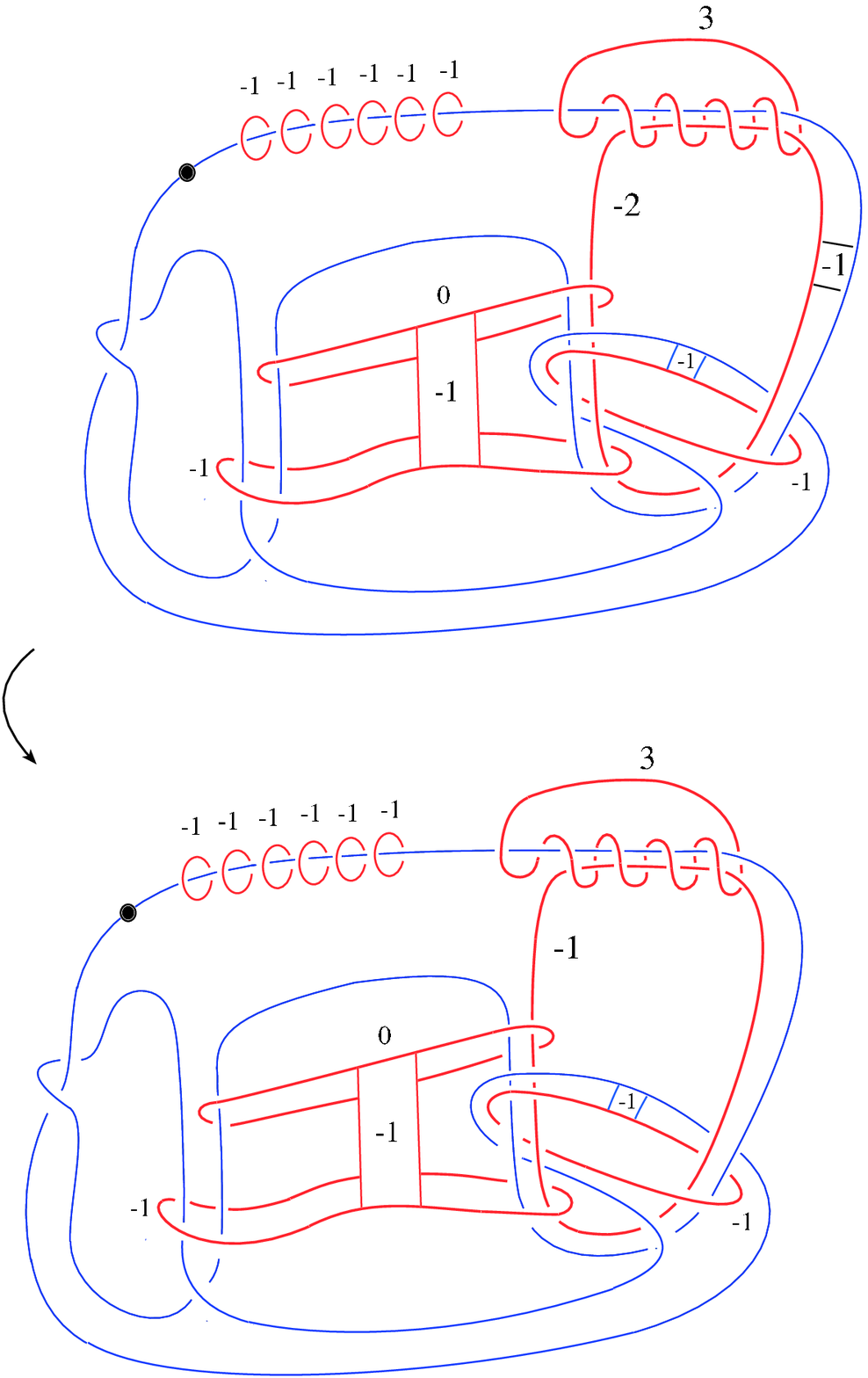}
\caption{}
\label{fig:poodles}
\end{figure}

\begin{figure}
\includegraphics[width=.70\textwidth]{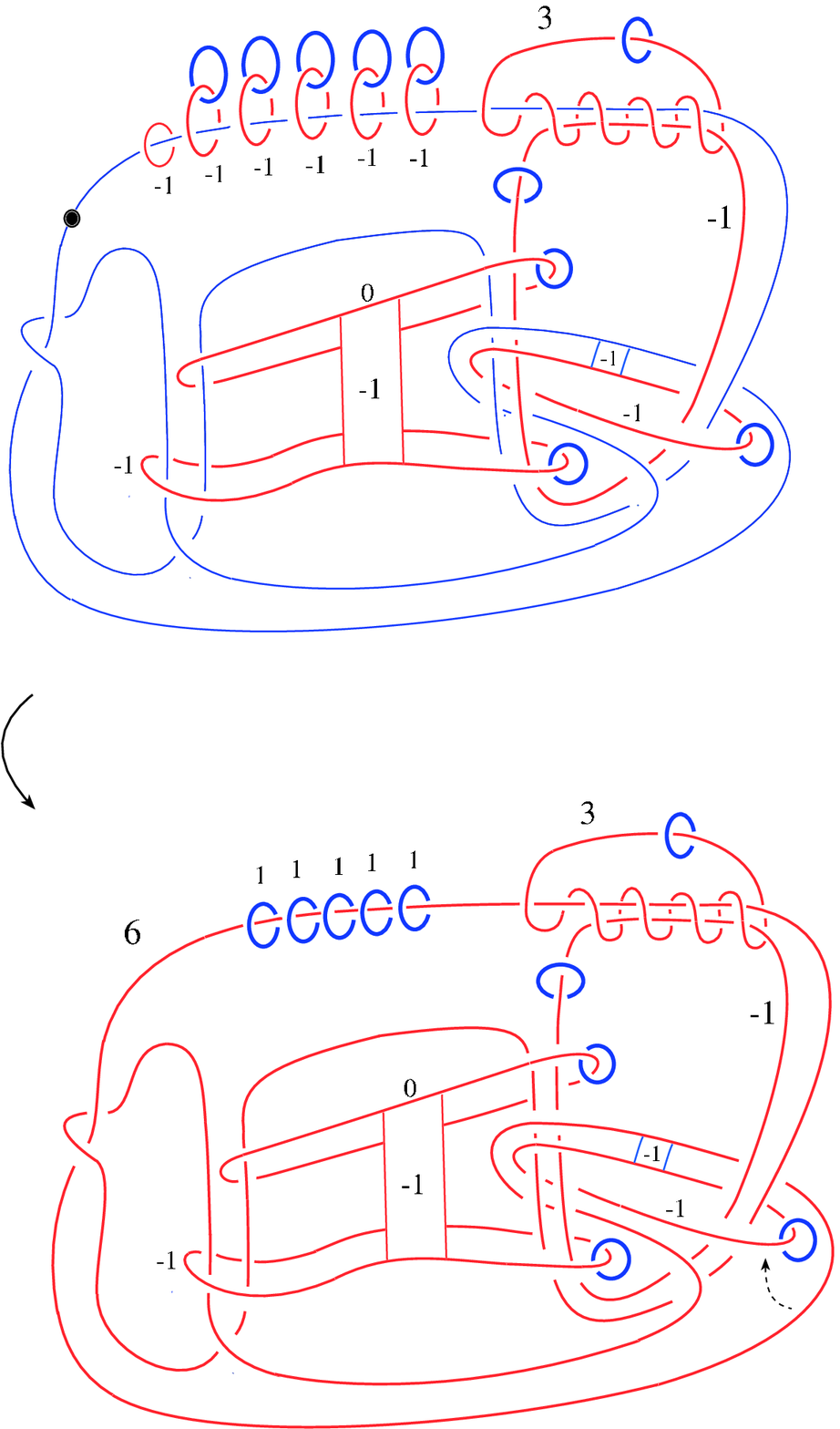}
\caption{}
\label{fig:poodles}
\end{figure}

\begin{figure}
\includegraphics[width=.80\textwidth]{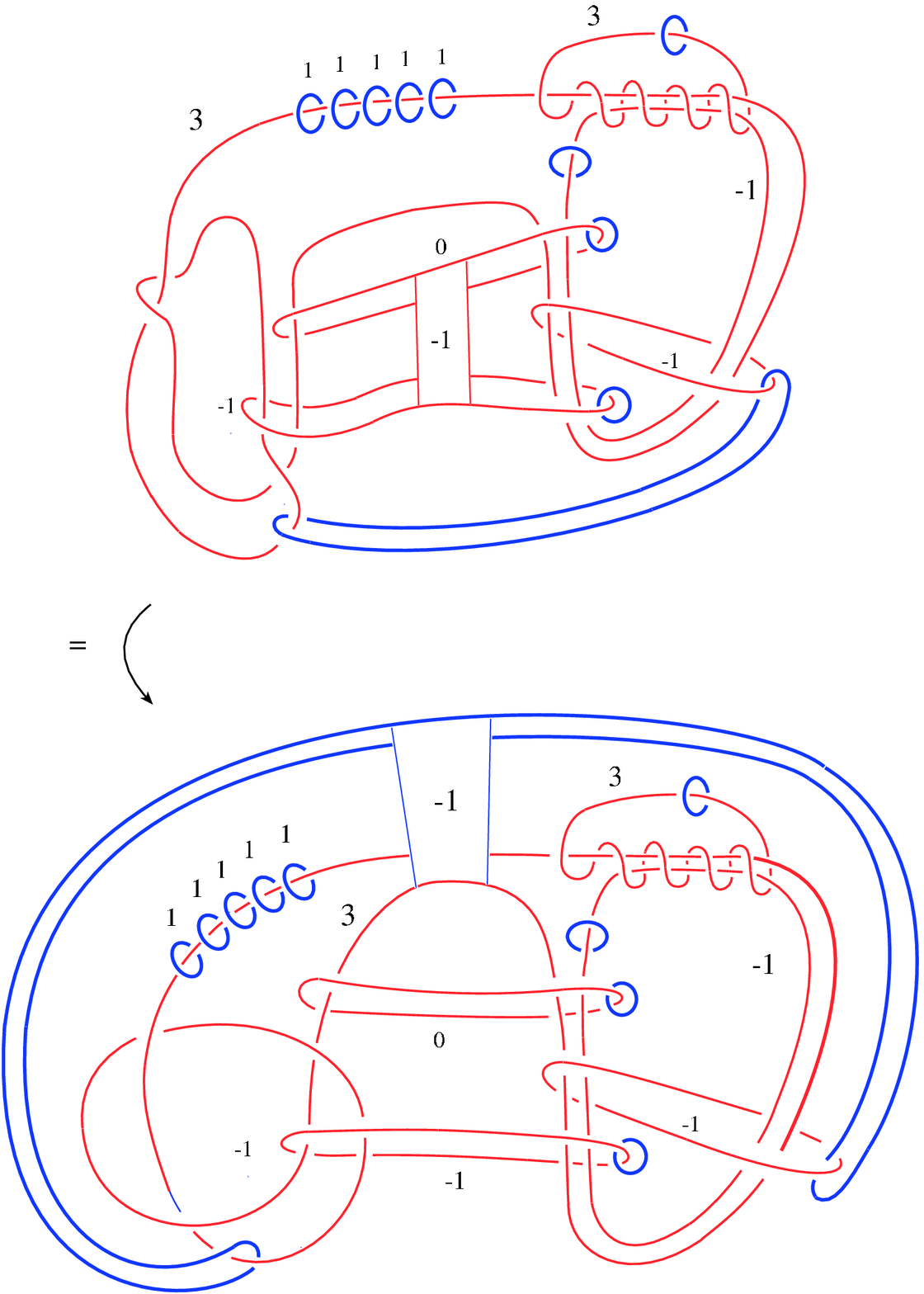}
\caption{}
\label{fig:poodles}
\end{figure}

\begin{figure}
\includegraphics[width=.75\textwidth]{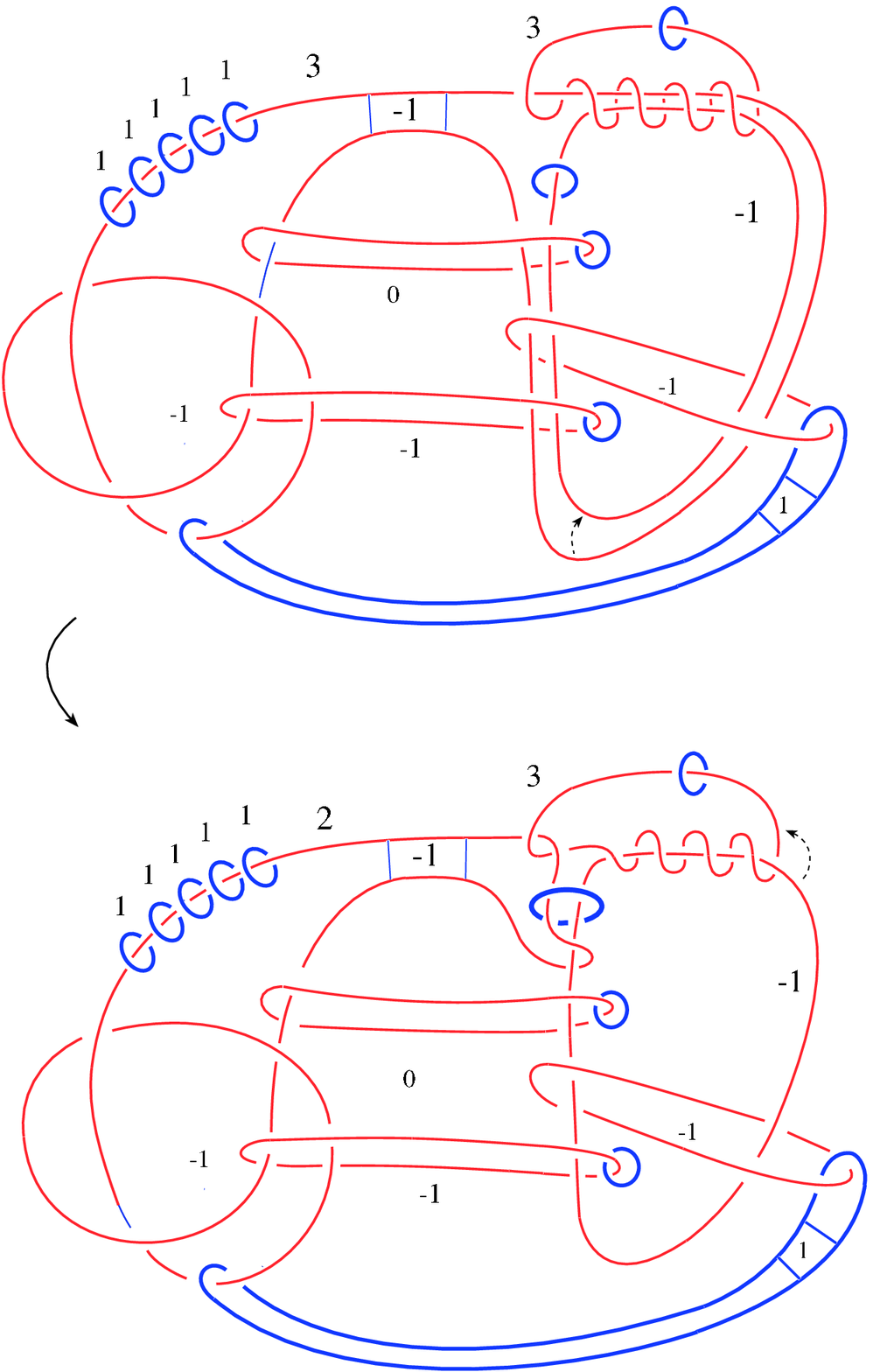}
\caption{}
\label{fig:poodles}
\end{figure}

\begin{figure}
\includegraphics[width=.75\textwidth]{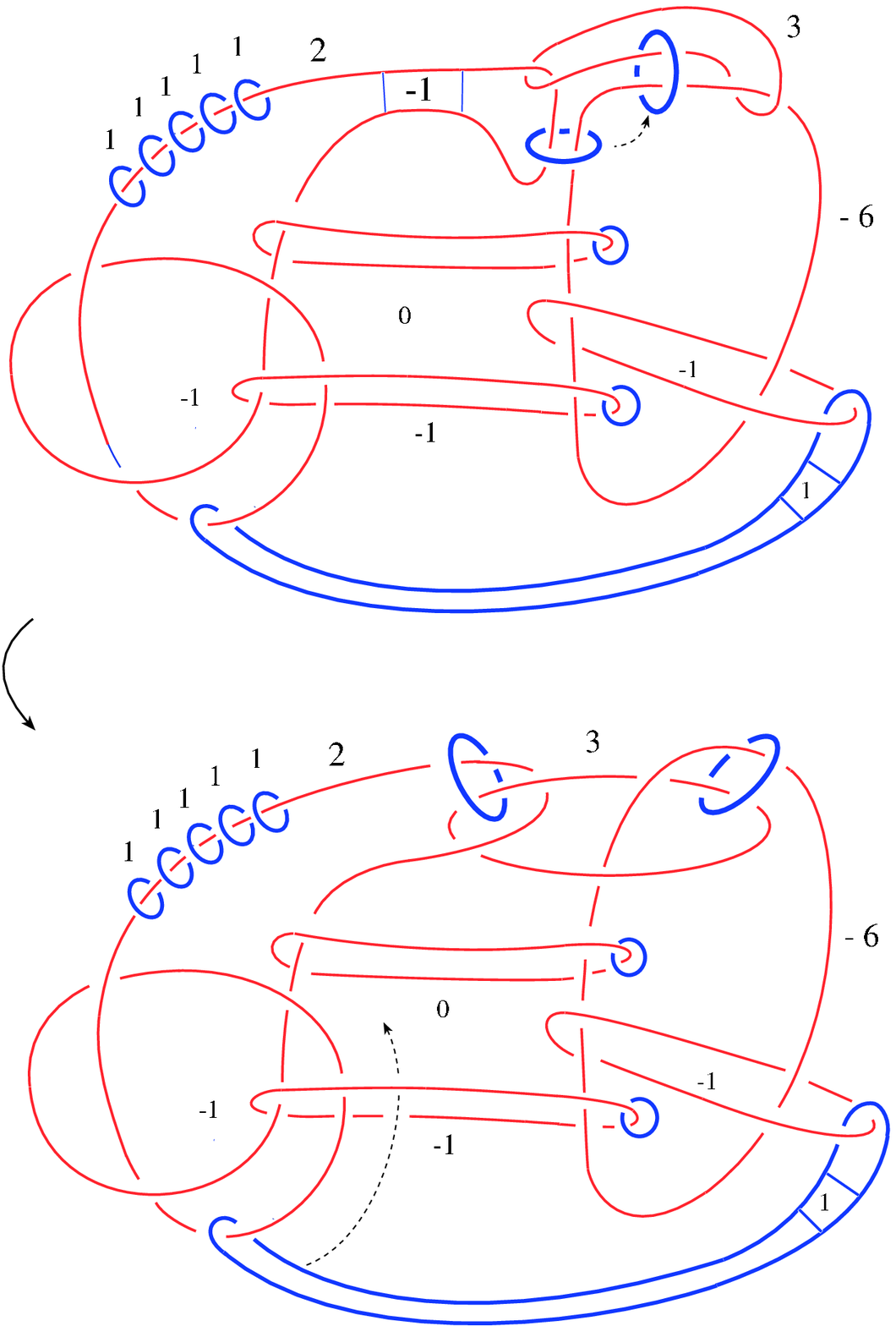}
\caption{}
\label{fig:poodles}
\end{figure}

\begin{figure}
\includegraphics[width=.75\textwidth]{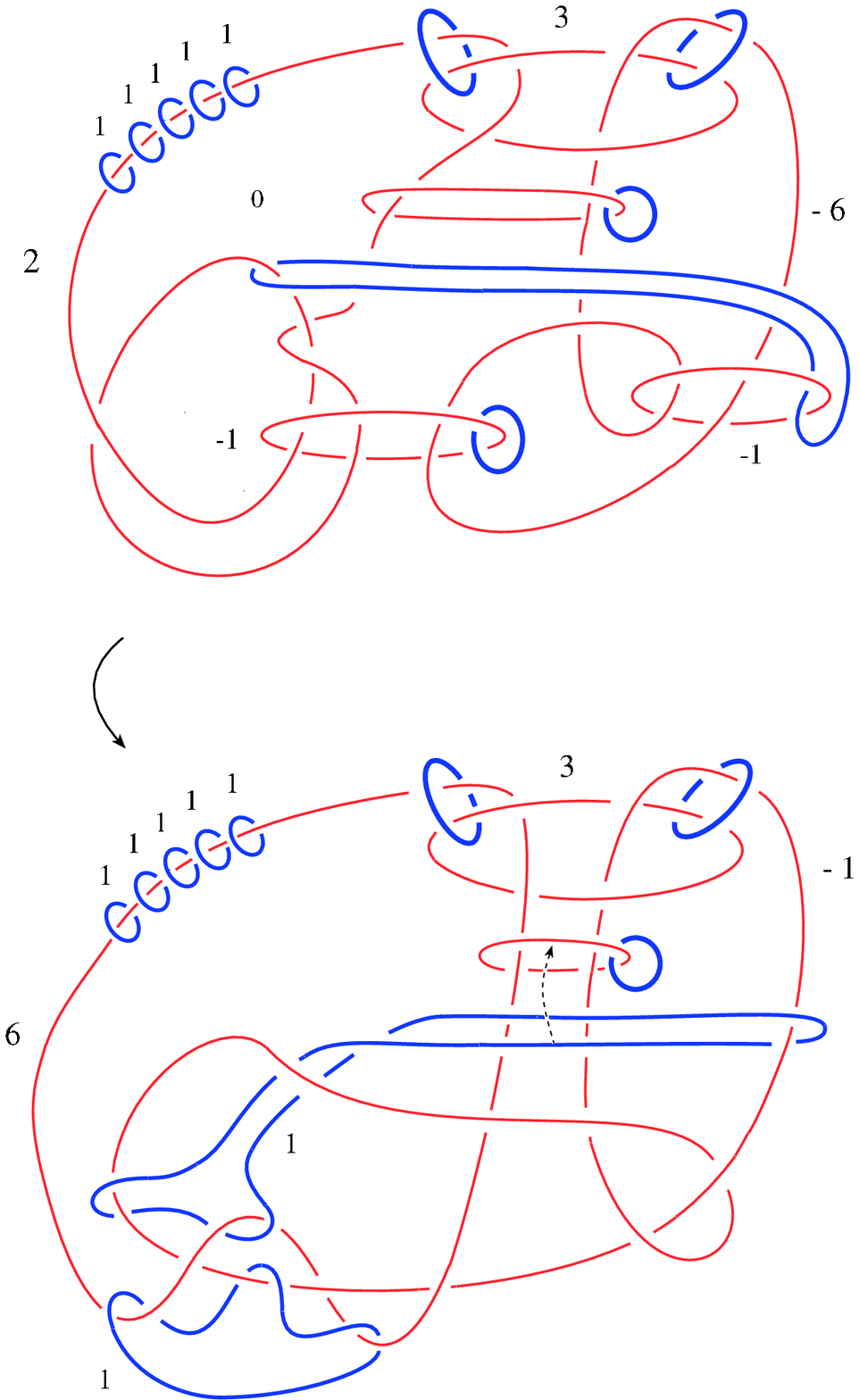}
\caption{}
\label{fig:poodles}
\end{figure}

\begin{figure}
\includegraphics[width=.75\textwidth]{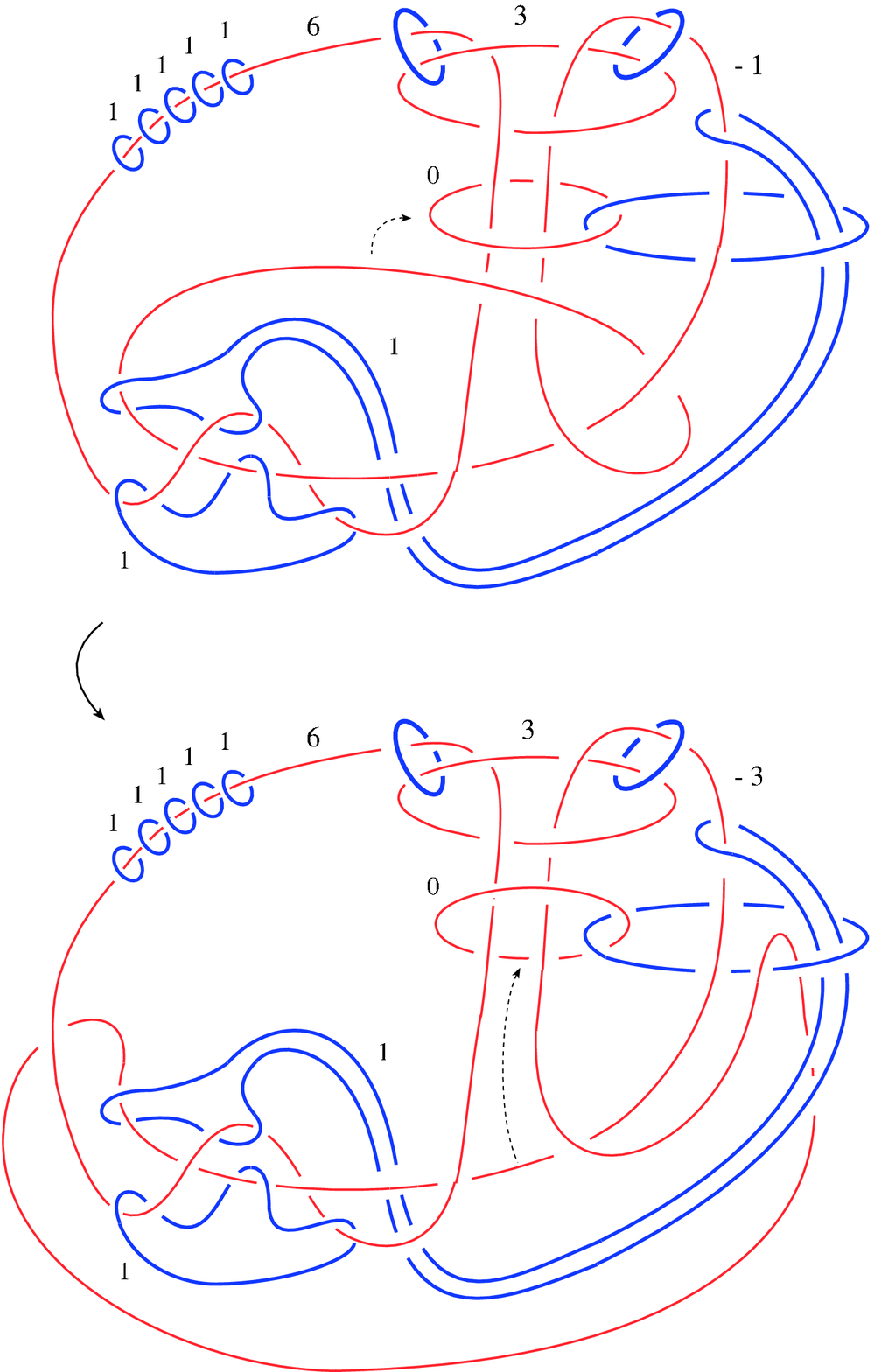}
\caption{}
\label{fig:poodles}
\end{figure}

\begin{figure}
\includegraphics[width=.75\textwidth]{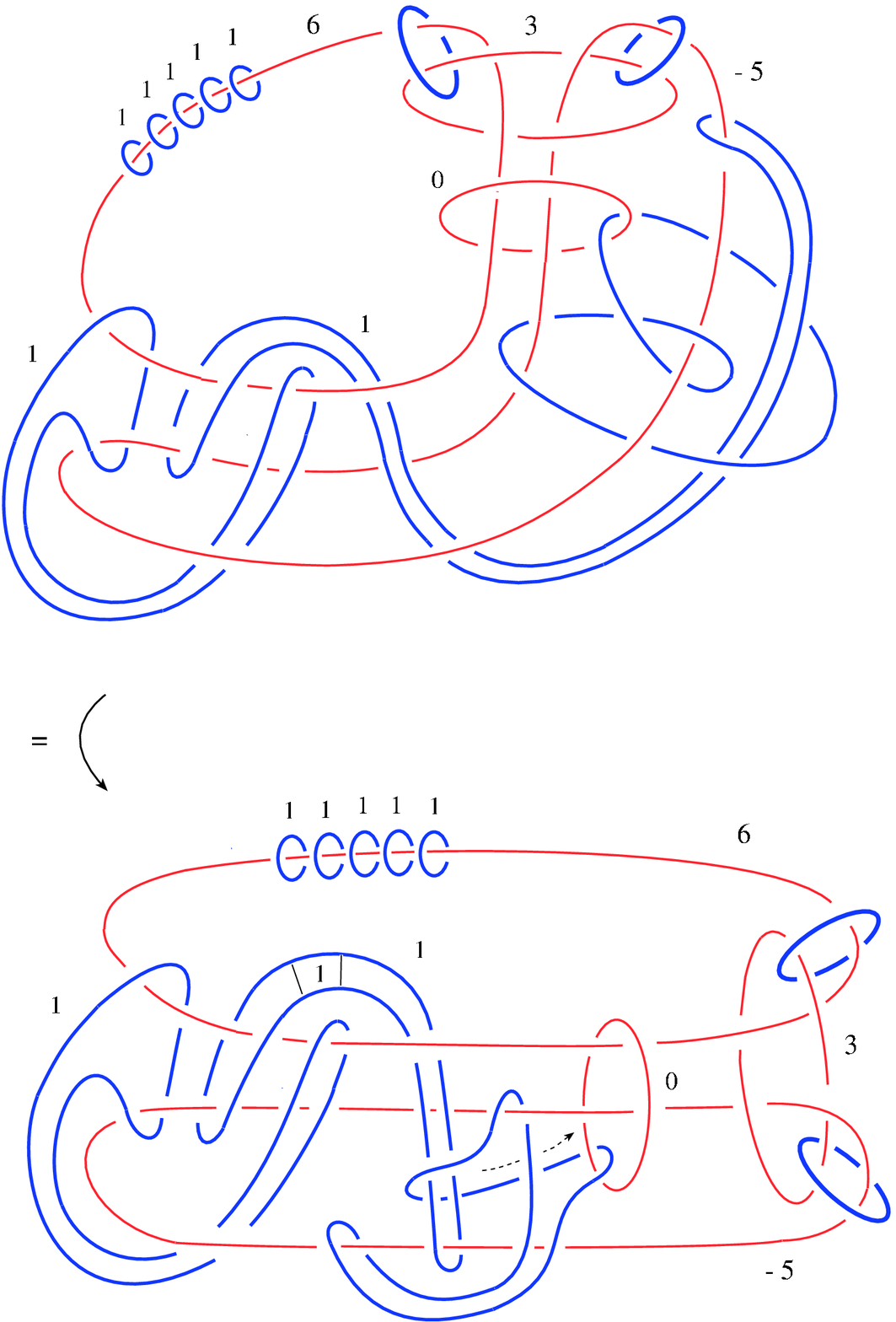}
\caption{}
\label{fig:poodles}
\end{figure}

\begin{figure}
\includegraphics[width=.85\textwidth]{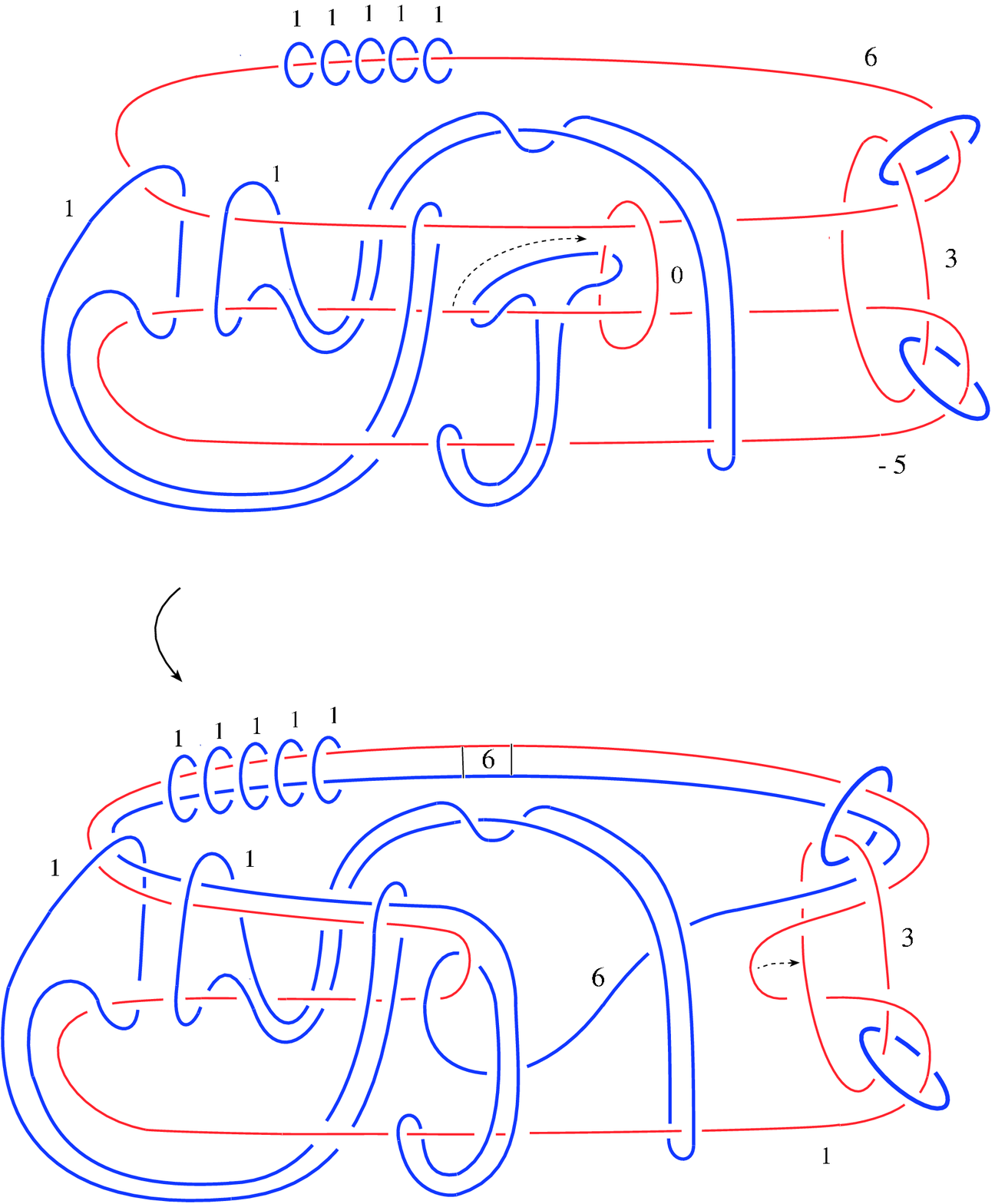}
\caption{}
\label{fig:poodles}
\end{figure}

\begin{figure}
\includegraphics[width=.75\textwidth]{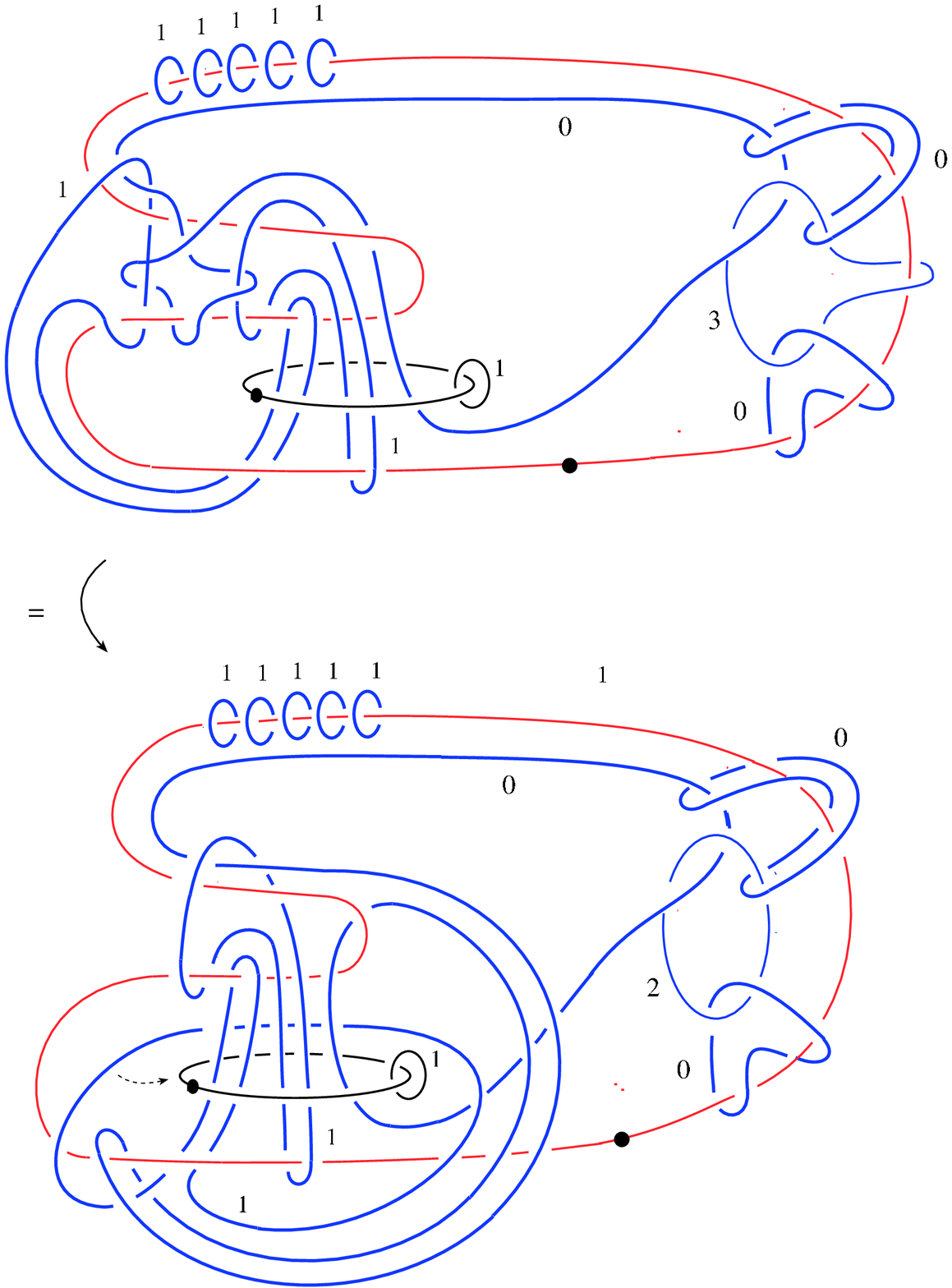}
\caption{}
\label{fig:poodles}
\end{figure}

\begin{figure}
\includegraphics[width=.75\textwidth]{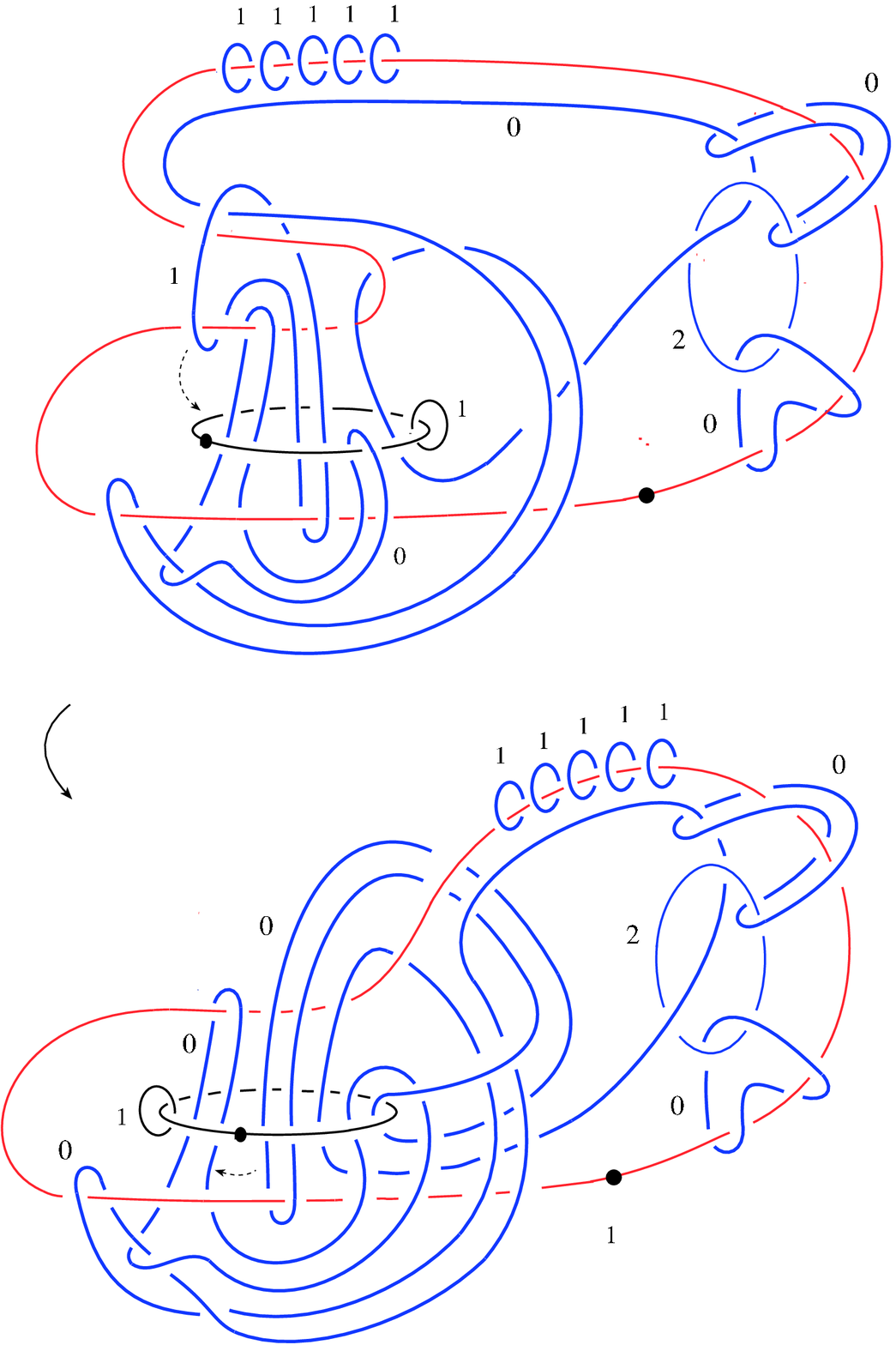}
\caption{}
\label{fig:poodles}
\end{figure}

\begin{figure}
\includegraphics[width=.85\textwidth]{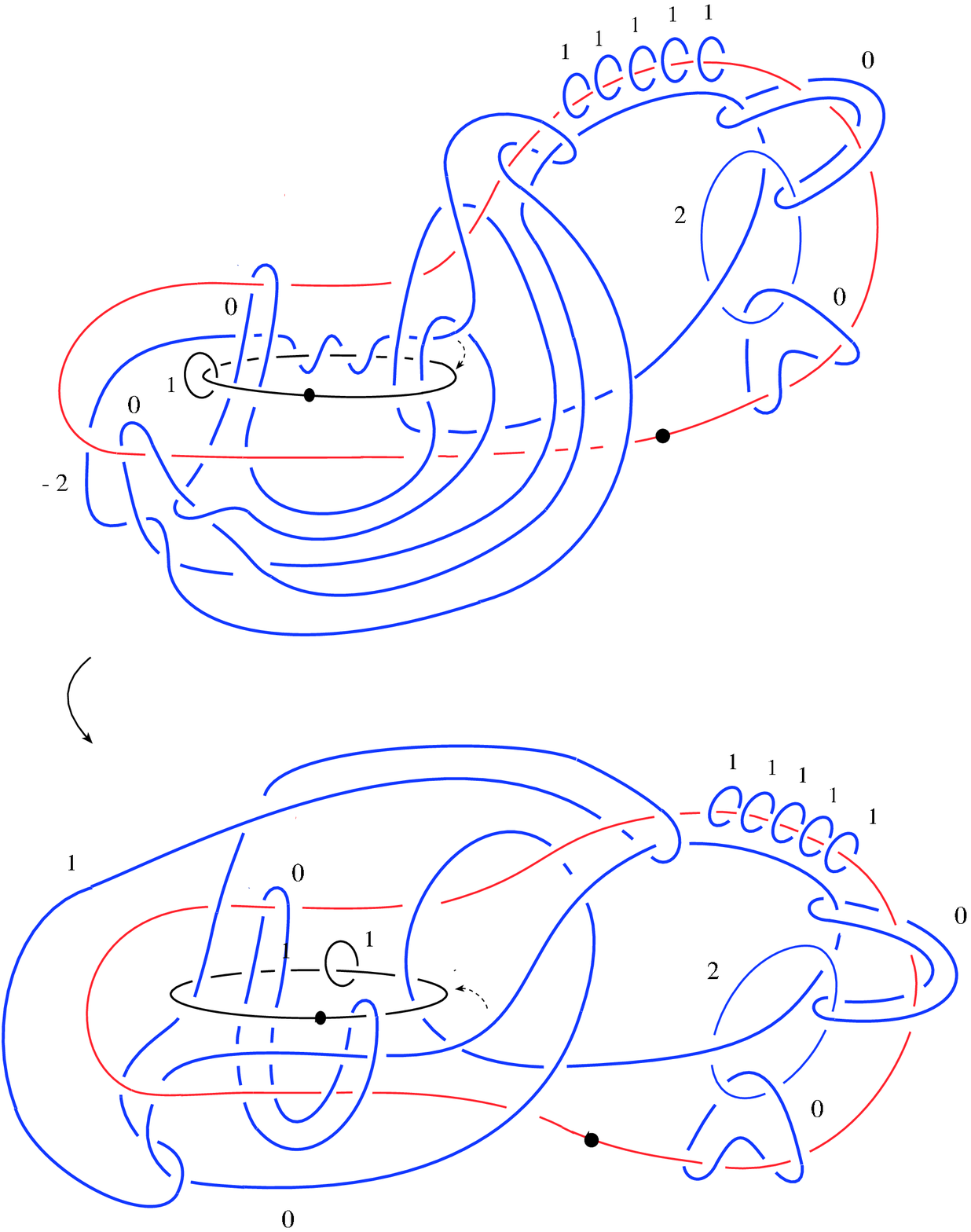}
\caption{}
\label{fig:poodles}
\end{figure}

\begin{figure}
\includegraphics[width=.75\textwidth]{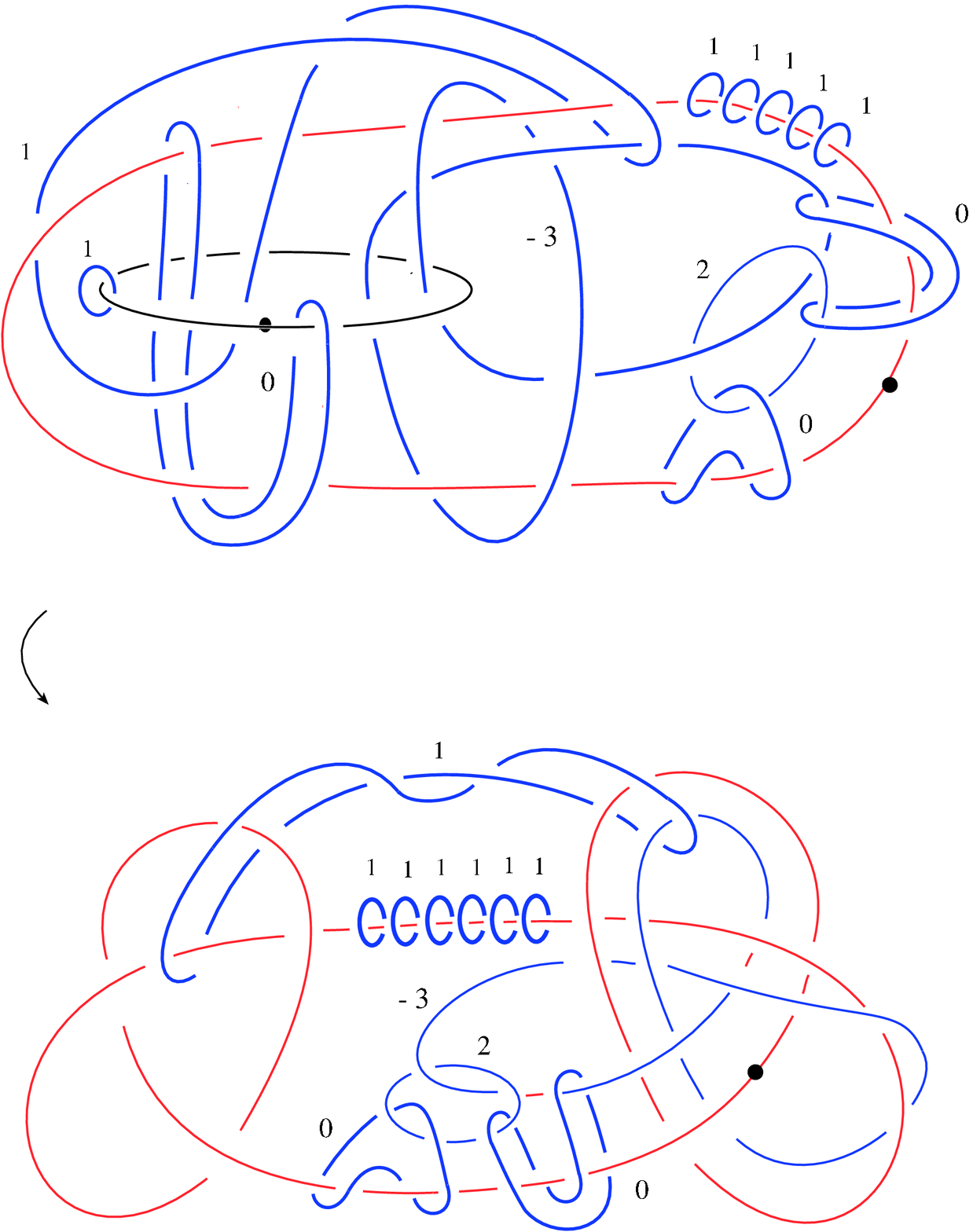}
\caption{}
\label{fig:poodles}
\end{figure}

\begin{figure}
\includegraphics[width=.75\textwidth]{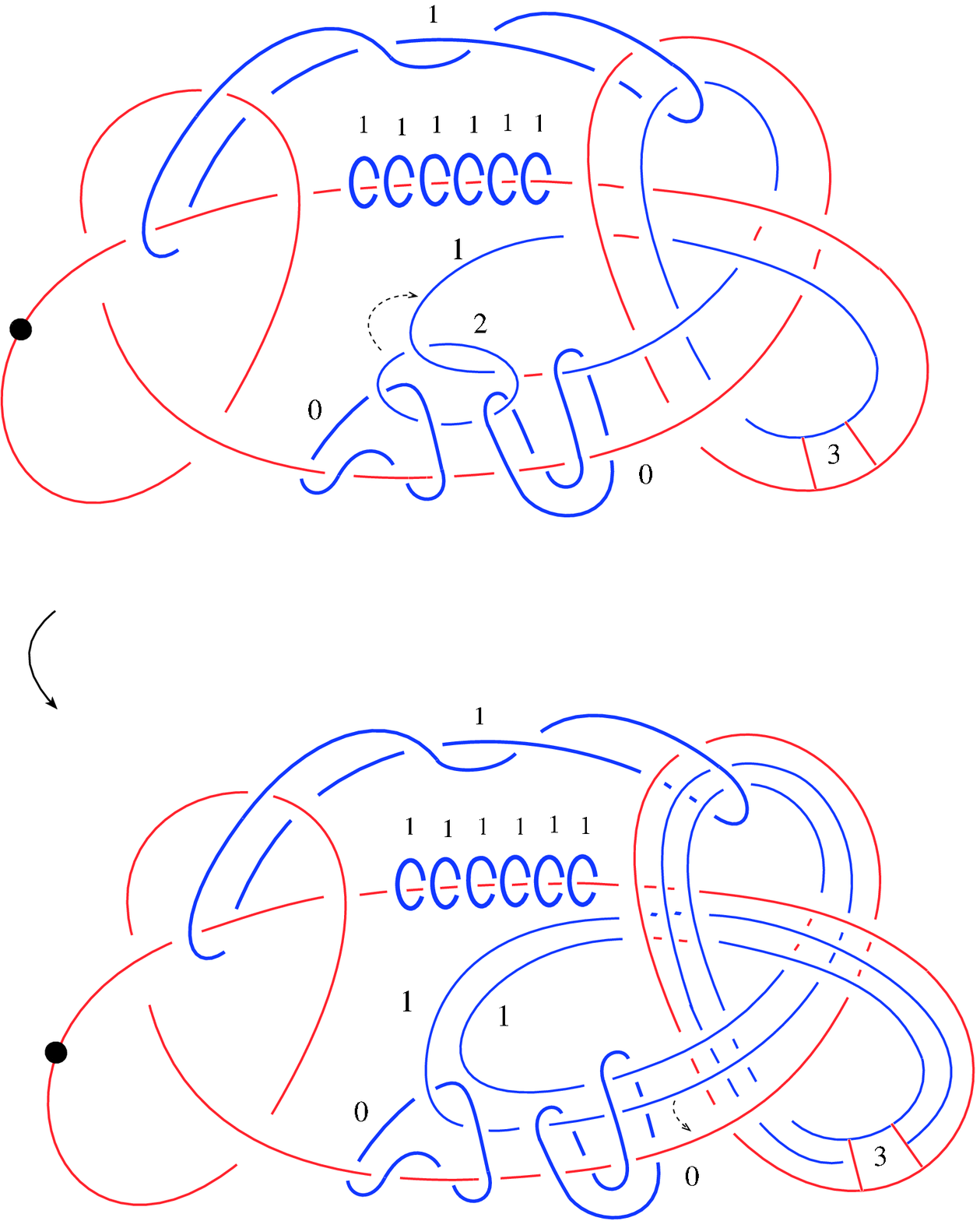}
\caption{}
\label{fig:poodles}
\end{figure}

\begin{figure}
\includegraphics[width=.75\textwidth]{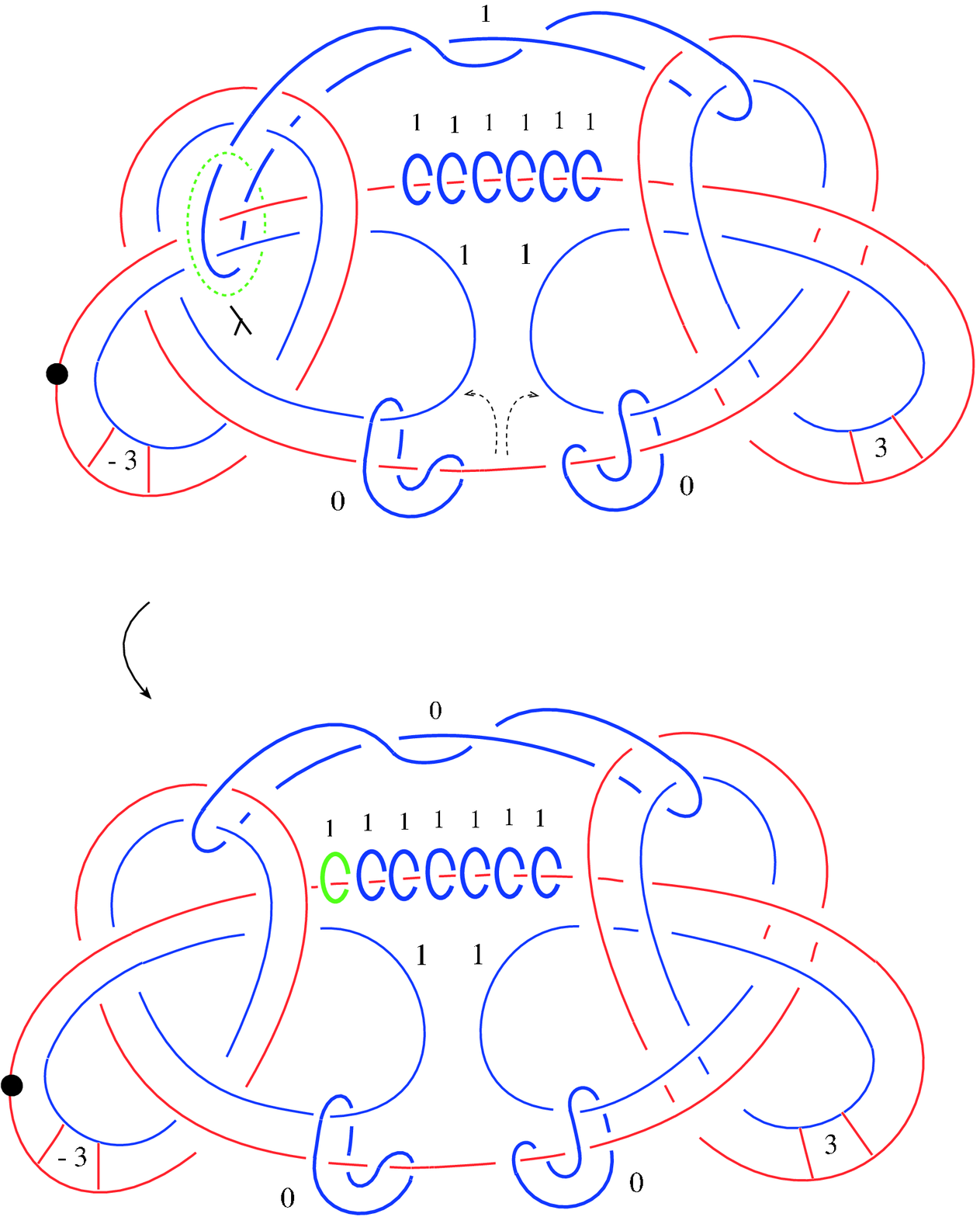}
\caption{}
\label{fig:poodles}
\end{figure}

\begin{figure}
\includegraphics[width=.70\textwidth]{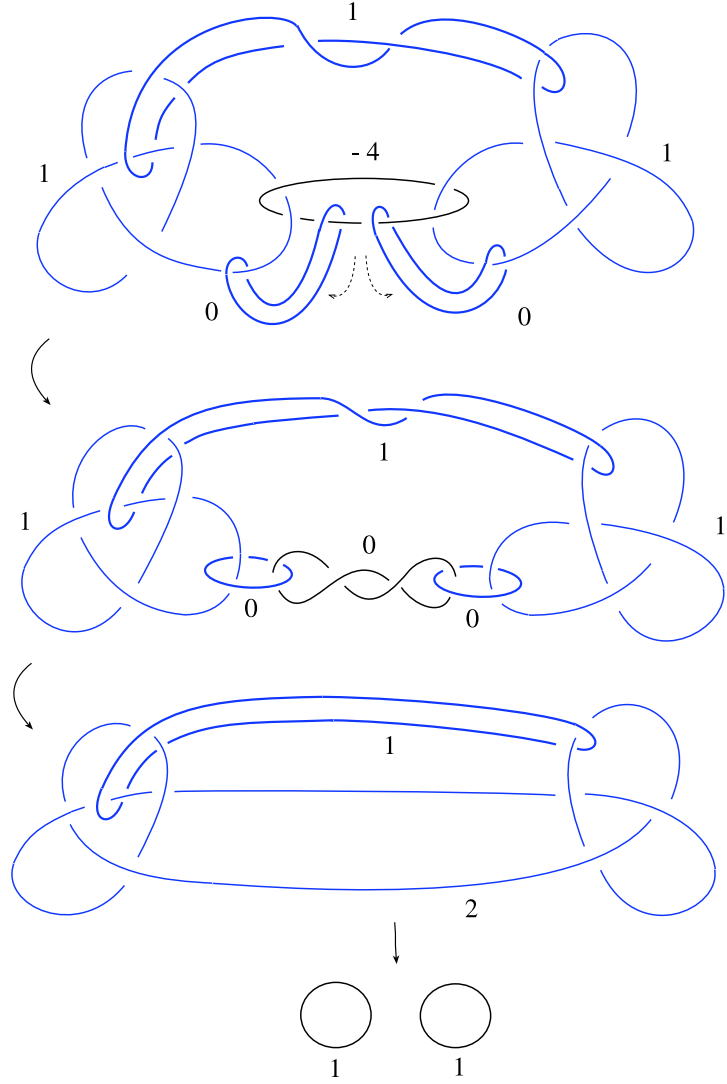}
\caption{The indicated handle slides  introduces a left twist below, canceling the right twist at the top handle.}
\label{fig:poodles}
\end{figure}

\begin{figure}
\includegraphics[width=.65\textwidth]{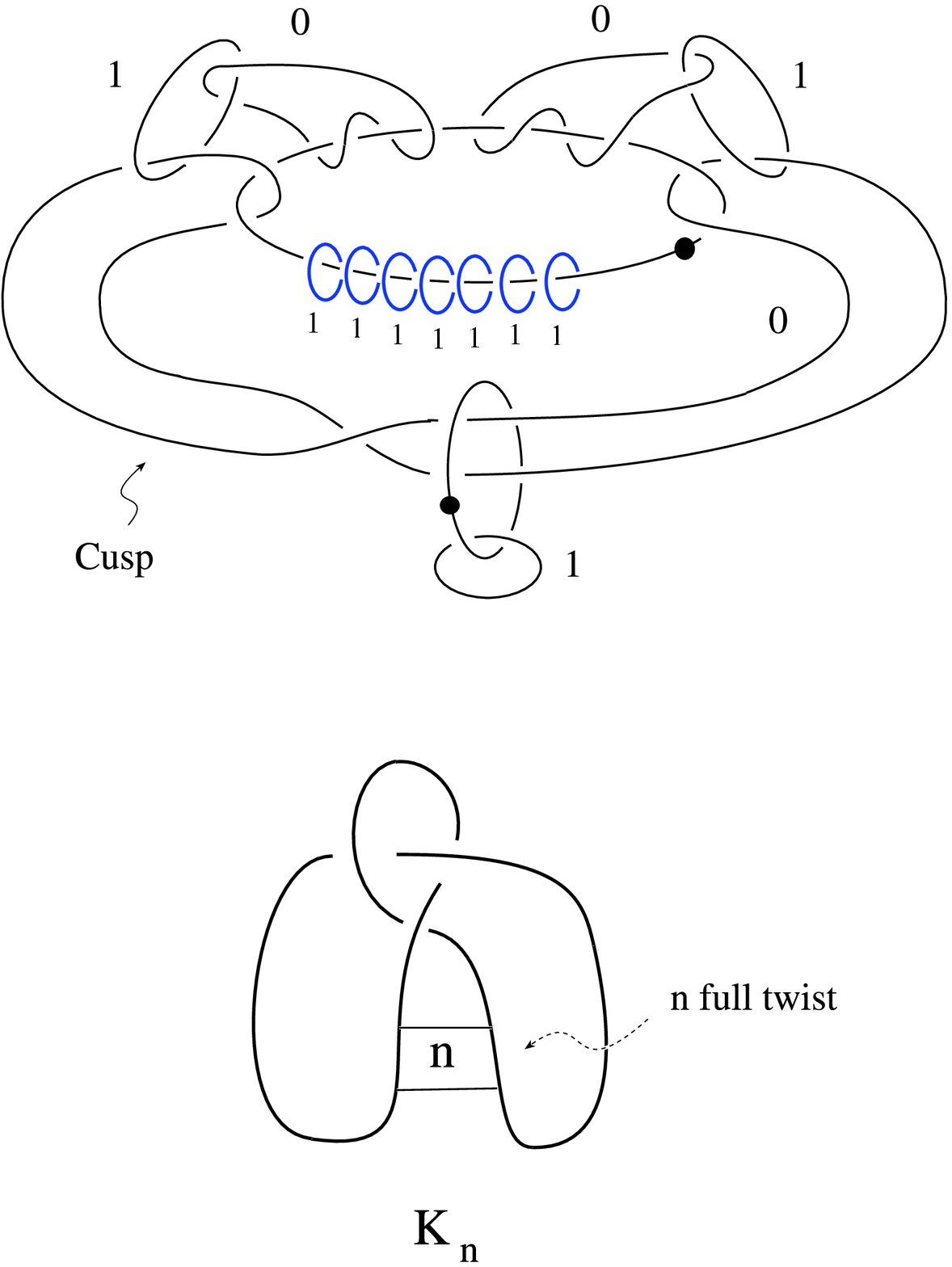}
\caption{}
\label{fig:poodles}
\end{figure}

\begin{figure}
\includegraphics[width=.75\textwidth]{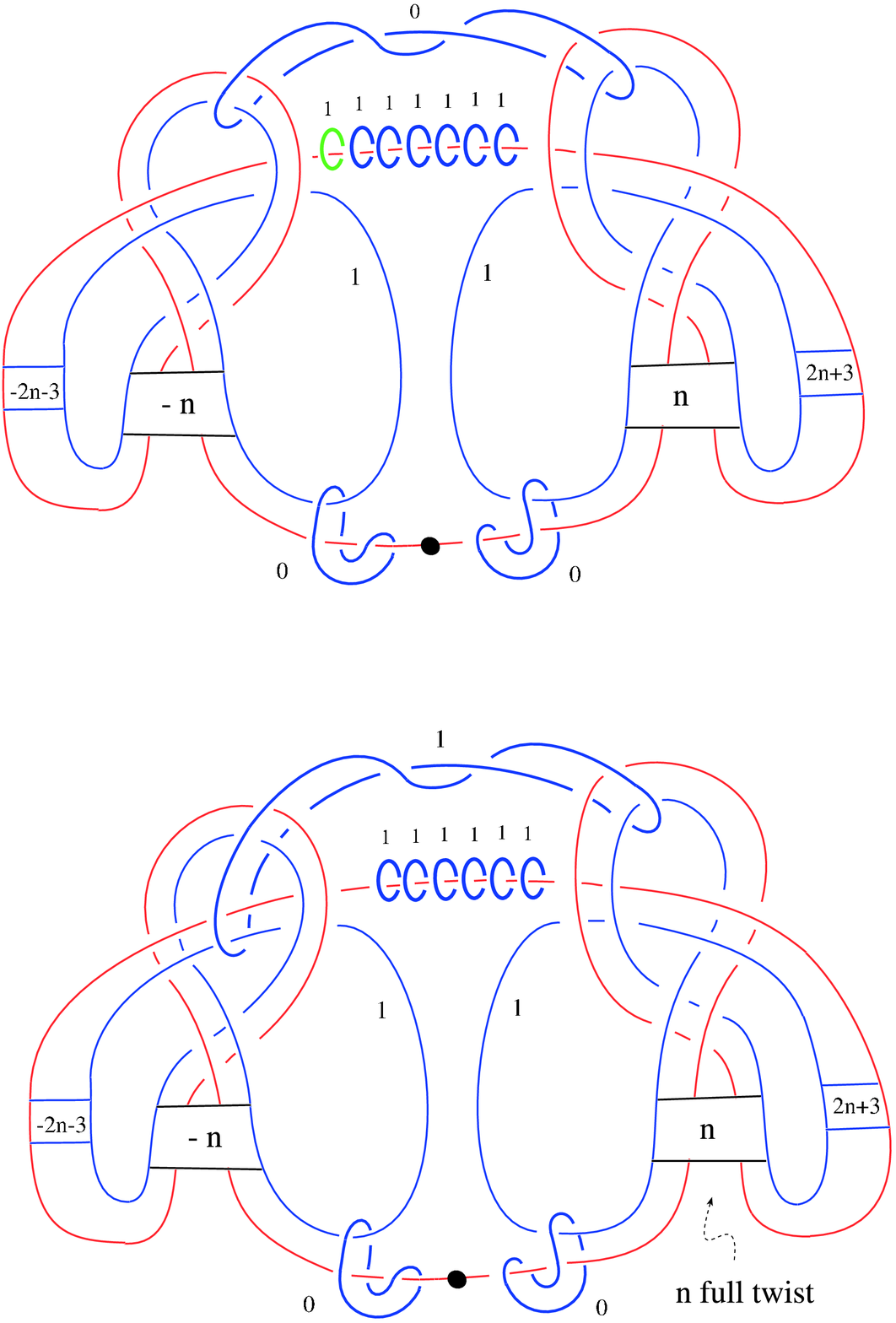} 
\caption{}
\label{fig:poodles}
\end{figure}

\end{document}